\documentclass[a4paper, 12pt]{amsart}
\usepackage{amsmath, amssymb, eucal, amscd, amstext, srctex}

\usepackage{enumerate}

\newtheorem{theorem}{Theorem}
\newtheorem{corollary}[theorem]{Corollary}
\newtheorem{lemma}[theorem]{Lemma}

\theoremstyle{definition}

\theoremstyle{remark}
\newtheorem{remark}[theorem]{Remark}

\newtheorem{example}[theorem]{Example}

\begin{document}

\title{Inner products and module maps of Hilbert $C^*$-modules}

\author{Ming-Hsiu Hsu \and Ngai-Ching Wong}

\address{Department of Applied Mathematics, National Sun Yat-sen University, Kaohsiung, 80424, Taiwan.}
\email[Ming-Hsiu Hsu]{hsumh@math.nsysu.edu.tw}
\email[Ngai-Ching Wong]{wong@math.nsysu.edu.tw}

\thanks{This work is jointly supported by a Taiwan NSC Grant ().}


\date{}


\begin{abstract}
Let $E$ and $F$ be two Hilbert $C^*$-modules over $C^*$-algebras $A$ and $B$, respectively.
Let $T$ be a surjective linear isometry from $E$ onto $F$ and $\varphi$ a map from $A$ into $B$.
We will prove in this paper that if the $C^*$-algebras $A$ and $B$ are commutative, then
$T$ preserves the inner products and $T$ is a module map, i.e., there exists a $*$-isomorphism $\varphi$ between
the $C^*$-algebras such that $$\langle Tx,Ty\rangle=\varphi(\langle x,y\rangle),$$ and
$$T(xa)=T(x)\varphi(a).$$
In case $A$ or $B$ is noncommutative $C^*$-algebra, $T$ may not satisfy the equations above in general.
We will also give some condition such that $T$ preserves the inner products and $T$ is a module map.
\end{abstract}

\maketitle

\section{Introduction}

A (right) Hilbert $C^*$-module over a $C^*$-algebra $A$ is a right $A$-module $E$ equipped
with $A$-valued inner product $\langle\cdot,\cdot\rangle$ which is conjugate $A$-linear in the first variable
and $A$-linear in the second variable such that $E$ is a Banach space with respect to the norm
$\|x\|=\|\langle x,x\rangle\|^{1/2}$.

Let $X$ be a locally compact Hausdorff space and $H$ a Hilbert space, the Banach space $C_0(X,H)$ of all
continuous $H$-valued functions vanishing at infinity is a Hilbert $C^*$-module over the $C^*$-algebra $C_0(X)$
with inner product $\langle f,g\rangle(x):=\langle f(x),g(x)\rangle$ and module operation $(f\phi)(x)=f(x)\phi(x)$,
for all $f\in C_0(X,H)$ and $\phi\in C_0(X)$. Every $C^*$-algebra $A$ is a Hilbert $C^*$-module
over itself with inner product $\langle a,b\rangle:=a^*b$.

Let $X$ and $Y$ be two locally compact Hausdorff spaces, the Banach-Stone theorem states that every surjective
linear isometry between $C_0(X)$ and $C_0(Y)$ is a weighted composition operator. More precisely, let $T$ be a
surjective linear isometry from $C_0(X)$ onto $C_0(Y)$, then there exists a continuous function $h\in C_0(Y)$
with $|h(y)|=1$, for all $y$ in $Y$, and a homeomorphism $\varphi$ from $Y$ onto $X$ such that $T$ is of the form:
\begin{eqnarray}\label{eq:wco}
Tf(y)=h(y)f(\varphi(y)), \forall f\in C_0(X), \forall y\in Y.
\end{eqnarray}
Let $H_1$ and $H_2$ be two Hilbert spaces. In \cite{Jerison50}, Jerison characterizes surjective linear isometries
between $C_0(X,H_1)$ and $C_0(Y,H_2)$, see also \cite{Lau,Jeang03}. It is said that every surjective linear isometry
$T$ from $C_0(X,H_1)$ onto $C_0(Y,H_2)$ is also of the form (\ref{eq:wco}) in which $h(y)$ is a unitary operator
from $H_1$ onto $H_2$ and $h$ is continuous from $Y$ into $(B(H_1,H_2),SOT)$, the space of all bounded linear
operators with the strong operator topology.
In this case, we can find a relationship of inner products of $C_0(X,H_1)$ and $C_0(Y,H_2)$ by a simple calculation:
\begin{eqnarray*}
\langle Tf, Tg\rangle(y)
&=&\langle Tf(y),Tg(y)\rangle=\langle h(y)(f(\varphi(y))), h(y)(f(\varphi(y)))\rangle\\
&=&\langle f(\varphi(y)), f(\varphi(y))\rangle=\langle f,g\rangle\circ\varphi(y).
\end{eqnarray*}
i.e. $$\langle Tf, Tg\rangle=\langle f,g\rangle\circ\varphi.$$
Let $R_{\varphi}:C_0(X)\rightarrow C_0(Y)$ be the $*$-isomorphism defined by $R_{\varphi}(\phi)=\phi\circ\varphi$.
Then $T$ preserves the inner products with respect to $R_{\varphi}$, i.e.,
$$\langle Tf, Tg\rangle=R_{\varphi}(\langle f,g\rangle).$$
By (\ref{eq:wco}), it is easy to see that $T$ is a module map with respect to $R_{\varphi}$ in the sense
$$T(f\phi)=T(f)R_{\varphi}(\phi),\ \mbox{for all $f\in C_0(X,H_1)$ and $\phi\in C_0(X)$}.$$
It is natural to ask if these properties are true for surjective linear isometries between Hilbert $C^*$-modules
over $C^*$-algebras. We will show in this paper that the answer is yes if the $C^*$-algebras are commutative.
Unfortunately, if one of the $C^*$-algebras is noncommutative, the answer is more complicated. We will give
an example (see Example \ref{ex:cex}) to explain this is not true in general. And we will give a
condition on $T$ (see Theorem \ref{thm:main}) such that $T$ is a module map and preserves the inner products.

\section{Preliminaries}
Let $E$ be a Hilbert $C^*$-module over $C^*$-algebra $A$. We set $\langle E,E\rangle$ to be the linear span of
elements of the form $\langle x,y\rangle$, $x,y\in E$.
$E$ is said to be \emph{full} if the closed two-sided ideal $\overline{\langle E,E\rangle}$ equal $A$.

A \emph{$JB^*$-triple} is a complex vector space $V$ with a
continuous mapping $V^3\rightarrow V,\ (x,y,z)\rightarrow\{x,y,z\}$,
called a \emph{Jordan triple product}, which is symmetric and linear
in $x, z$ and conjugate linear in $y$ such that for $x,y,z,u,v$ in
$V$, we have
\begin{enumerate}
    \item $\{x,y,\{z,u,v\}\}=\{\{x,y,z\},u,v\}-\{z,\{y,x,u\},v\}+\{z,u,\{x,y,v\}\}$;
    \item the mapping $z\rightarrow \{x,x,z\}$ is hermitian and has non-negative spectrum;
    \item $\|\{x,x,x\}\|=\|x\|^3$.
\end{enumerate}
In \cite{Isidro03}, J. M. Isidro shows that every Hilbert $C^*$-module is
a ${\rm JB}^*$-triple with the Jordan triple product
$$
\{x,y,z\}=\displaystyle\frac{1}{2}(x\langle y,z\rangle+z\langle y,x\rangle).
$$
A well-known theorem of Kaup \cite{Kaup83} (see also \cite{Chu05}) states that every surjective linear isometry
between ${\rm JB}^*$-triples is a \emph{Jordan triple homomorphism}, i.e., it preserves the Jordan triple product
$$T\{x,y,z\}=\{Tx,Ty,Tz\}, \forall x,y,z\in E.$$
Hence, if $T$ is a surjective linear isometry between Hilbert $C^*$-modules, then
\begin{eqnarray}\label{eq:triple}
T(x\langle y,z\rangle+z\langle y,x\rangle)=Tx\langle Ty,Tz\rangle+Tz\langle Ty,Tx\rangle, \forall x,y,z\in E.
\end{eqnarray}
The equation (\ref{eq:triple}) holds if and only if
\begin{eqnarray}\label{eq:tri}
T(x\langle x,x\rangle)=Tx\langle Tx,Tx\rangle, \forall x\in E
\end{eqnarray}
by triple polarization
$$2\{x,y,z\}=\frac{1}{8}\sum\limits_{\alpha^4=\beta^2=1}
\alpha\beta \langle x+\alpha y+\beta z,x+\alpha y+\beta z\rangle(x+\alpha y+\beta z).$$

A ternary ring of operators (TRO) between two Hilbert spaces $H$ and $K$ is a linear subspace $\mathfrak{R}$
of $B(H,K)$, the space of all bounded linear operators from $H$ into $K$, satisfying $AB^*C\in\mathfrak{R}$.
Zettl shows in \cite{Zettl83} that every Hilbert $C^*$-module is isomorphic to a norm closed TRO.
In this case, Hilbert $C^*$-modules have another \emph{triple product}, i.e.,
$$\{x,y,z\}:=x\langle y,z\rangle.$$
A map between TROs is said to be a \emph{triple homomorphism} if it preserves the triple products.
In the case of Hilbert $C^*$-modules, a map $T$ is a triple homomorphism if it satisfies
\begin{eqnarray}\label{eq:tri}
T(x\langle y,z\rangle)=Tx\langle Ty,Tz\rangle, \forall x,y,z.
\end{eqnarray}
We have known every surjective linear isometry
is a Jordan triple homomorphism, but it could not be a \emph{triple homomorphism},
see Example \ref{ex:cex}.

Let $\mathcal{R}$ be a TRO. Then $M_n(\mathcal{R})$, the space of all $n\times n$ matrices whose entries are
in $\mathcal{R}$, has a TRO-structure. Let $T$ be a map between TROs $\mathcal{R}_1$ and $\mathcal{R}_2$.
For all positive integer $n$, define maps $T^{(n)}:M_n(\mathcal{R}_1)\rightarrow M_n(\mathcal{R}_2)$ by
$T^{(n)}((x_{ij})_{ij})=(T(x_{ij}))_{ij}$. We call $T$ $n$-isometry if $T^{(n)}$ is isometric and
\emph{complete isometry} if each $T^{(n)}$ is isometric for all $n$.
It has been shown that a surjective linear isometry between TROs is a triple homomorphism if and only if
it is completely isometric. More details about TROs mentioned above, we refer to \cite{Zettl83},
see also \cite{NealRusso03, Hamana99}. In fact, Solel shows in \cite{Solel01} that every surjective 2-isometry
between two full Hilbert $C^*$-modules is necessarily completely isometric.

\section{Results}
Note that in the case of a commutative $C^*$-algebra $A=C_0(X)$, for some locally compact Hausdorff space $X$,
Hilbert $C^*$-modules over $C_0(X)$ are the same as Hilbert bundles, or equivalently, continuous fields of
Hilbert spaces, over $X$.

We showed the following theorem in \cite{HsuWong11}.

\begin{theorem}\label{thm:commu}
Let $E$ and $F$ be two Hilbert $C^*$-modules over commutative $C^*$-algebras $C_0(X)$ and $C_0(Y)$, respectively.
Then every surjective linear isometry from $E$ onto $F$ is a weighted composition operator
$$Tf(y)=h(y)(f(\varphi(y))), \forall f\in E, \forall y\in Y$$
Here, $\varphi$ is a homeomorphism from $Y$ onto $X$, $h(y)$ is a unitary operator between the corresponding fibers
of $E$ and $F$, for all $y$ in $Y$.
\end{theorem}

By the similar argument discussed in the introduction, we have

\begin{corollary}
Every surjective linear isometry between Hilbert $C^*$-modules over commutative $C^*$-algebras preserves the inner
products and is a module map.
\end{corollary}

Now we discuss the case of noncommutative $C^*$-algebras.
From equation (\ref{eq:tri}), it seems that a surjective linear isometry $T$ indicates that $T$ preserves
inner products and that $T$ is a module map. We explain this could be not true in general by a example.

\begin{example}\label{ex:cex}
Given a positive integer $n$.
The Hilbert column space $H_c^n$ is the subspace of $M_n(\mathbb{C})$ consisting of all matrices whose non-zero entries
are only in the first column. Similarly, the Hilbert row space is the subspace consisting of matrices whose non-zero
entries are only in the first row. Clearly, $H_c$ and $H_r$ are right Hilbert $C^*$-modules over $C^*$-algebras
$\mathbb{C}$ and $M_n(\mathbb{C})$, respectively, with the inner product $\langle A,B\rangle:=A^*B$.
Define a surjective linear isometry $T:H_r^n\rightarrow H_c^n$ by $T(A)=A^t$, the transpose of $A$.
Then $\langle T(A), T(B)\rangle=tr\langle A, B\rangle$, the trace of $\langle A, B\rangle$,
but $T$ is not a module map with respect to the trace.
For the surjective linear isometry $T:H_c^n\rightarrow H_r^n$, $T(A)=A^t$.
Let $\varphi:\mathbb{C}\rightarrow M_n(\mathbb{C})$ be defined by $\varphi(\alpha)=\alpha I$.
Then $T$ is a module map with respect to $\varphi$, but the equation $\langle TA,TB\rangle=\varphi(\langle A,B\rangle)$
does not hold. It is clear that $T$ does not satisfy the equation (\ref{eq:tri}).
\end{example}

\begin{remark}
In fact, the corollary above says that there exists a $*$-isomorphism $\varphi$ between the $C^*$-algebras such that
$$\langle Tx,Ty\rangle=\varphi(\langle x,y\rangle)$$ and $$T(xa)=T(x)\varphi(a).$$ We have seen in the
Example \ref{ex:cex} that even if $T$ is a module map or preserves the inner products, the map $\varphi$ might be
just a linear map.
\end{remark}

In the following, $E$ and $F$ stand for two Hilbert $C^*$-modules over $C^*$-algebras $A$ and $B$, respectively.
$T$ is a map from $E$ into $F$ and $\varphi$ is a map from $A$ into $B$.
The following lemmas explain the relations of $T$, $\varphi$, when $T$ preserves the inner products and
when $T$ is a module map, see also \cite{Joita07}.

\begin{lemma}
If $\varphi$ is linear, every map $T$ from $E$ into $F$ which preserves the inner products with respect to $\varphi$
is linear.
\end{lemma}
\proof
Since $T$ preserves the inner products with respect to $\varphi$. Then for all $x, y$ and $z$ in $E$,
$\alpha$ in $\mathbb{C}$,
$$\langle T(\alpha x+y),Tz\rangle=\varphi(\langle \alpha x+y,z\rangle)
=\alpha\varphi(\langle x,z\rangle)+\varphi(\langle y,z\rangle)=\langle \alpha Tx+Ty,Tz\rangle.$$
Similarly, we have $$\langle Tx, T(\alpha y+z)\rangle=\langle Tx,\alpha Ty+Tz\rangle.$$
It is easy to show that
$$\langle T(\alpha x+y)-(\alpha Tx+Ty),T(\alpha x+y)-(\alpha Tx+Ty)\rangle=0.$$
This proves $T(\alpha x+y)=\alpha Tx+Ty$ and hence $T$ is linear.
\endproof

\begin{lemma}[\cite{Joita07}]
Let $T$ be a surjective linear map which preserves the inner products and is a module map w.r.t. $\varphi$.
If $F$ is full, then $\varphi$ is a $*$-homomorphism.
\end{lemma}
\proof
Let $a_1, a_2$ in $A$ and $\alpha$ in $\mathbb{C}$. It is easy to show that
\begin{eqnarray*}
& &T(x)(\varphi(\alpha a_1+a_2)-\alpha\varphi(a_1)-\varphi(a_2))\\
&=&T(x)\varphi(\alpha a_1+a_2)-\alpha T(x)\varphi(a_1)-T(x)\varphi(a_2)\\
&=&T(\alpha xa_1+xa_2)-\alpha T(xa_1)-T(xa_2)=0.
\end{eqnarray*}
and
\begin{eqnarray*}
& &T(x)(\varphi(a_1a_2)-\varphi(a_1)\varphi(a_2))\\
&=&T(x)\varphi(a_1a_2)-T(x)\varphi(a_1)\varphi(a_2)\\
&=&T(xa_1a_2)-T(xa_1a_2)=0.
\end{eqnarray*}
Since $T$ is surjective and $F$ is full, we have $\varphi(\alpha a_1+a_2)=\alpha\varphi(a_1)+\varphi(a_2)$
and $\varphi(a_1a_2)=\varphi(a_1)\varphi(a_2)$.

Let $x, y$ in $A$, we have
$$\varphi(\langle x,y\rangle^*)=\varphi(\langle y,x\rangle)=\langle Ty,Tx\rangle=\langle Tx,Ty\rangle^*
=\varphi(\langle x,y\rangle)^*.$$
For $a$ in $A$,
\begin{eqnarray*}
& &\langle T(x)(\varphi(a^*)-\varphi(a)^*),T(x)(\varphi(a^*)-\varphi(a)^*)\rangle\\
&=&\varphi(a^*)^*\varphi(\langle x,x\rangle)\varphi(a^*)-\varphi(a^*)^*\varphi(\langle x,x\rangle)\varphi(a)^*
-\varphi(a)\varphi(\langle x,x\rangle)\varphi(a^*)+\varphi(a)\varphi(\langle x,x\rangle)\varphi(a)^*\\
&=&(\varphi(\langle xa^*,x\rangle)\varphi(a^*))^*-(\varphi(a)\varphi(\langle x,x\rangle)\varphi(a^*))^*
-\varphi(\langle xa^*,xa^*\rangle)+(\varphi(a)\varphi(\langle x,xa^*\rangle))^*\\
&=&0.
\end{eqnarray*}
Hence, $T(x)(\varphi(a^*)-\varphi(a)^*)=0$ for all $x$ in $E$. Since $T$ is surjective and $F$ is full, we have
$\varphi(a^*)=\varphi(a)^*$.
\endproof

\begin{lemma}\label{lem:mp}
If $\varphi$ is a $*$-homomorphism, then every map $T$ which preserves the inner products w.r.t. $\varphi$ is a
module map w.r.t. $\varphi$.
\end{lemma}
\proof
Let $x$ and $y$ in $E$ and $a$ in $A$. Then
$$\langle T(xa),Ty\rangle=\varphi(\langle xa,y\rangle)=\varphi(a)^*\varphi(\langle x,y\rangle)
=\langle T(x)\varphi(a),Ty\rangle.$$
Similarly, we have $$\langle T(x),T(ya)\rangle=\langle T(x),T(y)\varphi(a)\rangle.$$
It is easy to show that
$$\langle T(xa)-T(x)\varphi(a),T(xa)-T(x)\varphi(a)\rangle=0.$$
Hence, $T(xa)=T(x)\varphi(a)$.
\endproof

\begin{lemma}[\cite{MS}]
Let $T$ be a surjective linear isometry and $\varphi$ a $*$-isomorphism.
If $T$ is a module map w.r.t. $\varphi$, then $T$ preserves the inner products with respect to $\varphi$.
\end{lemma}
\proof
It suffices to prove that $\langle Tx,Tx\rangle=\varphi(\langle x,x\rangle)$ for all $x$ in $E$.
Note that $|a|:=(a^*a)^{1/2}$. For all $b$ in $B$, let $\varphi(a)=b$, then
\begin{eqnarray*}
& &\||Tx|b\|^2=\|b^*|Tx|^2b\|=\|\langle T(x)\varphi(a),T(x)\varphi(a)\rangle\|\\
&=&\|\langle T(xa),T(xa)\rangle\|=\|\langle xa,xa\rangle\|=\||x|a\|^2=\|\varphi(|x|a)\|^2=\|\varphi(|x|)b\|^2.
\end{eqnarray*}
By Lemma 3.5 in \cite{Lance}, we get $|Tx|=(\varphi(|x|)$ and
hence $\langle Tx,Tx\rangle=\varphi(\langle x,x\rangle).$
\endproof

\begin{theorem}\label{thm:main}
Let $T$ be a surjective linear 2-isometry from $E$ onto $F$. Then there exists a $*$-isomorphism $\varphi$ from
$\overline{\langle E,E\rangle}$ onto $\overline{\langle F,F\rangle}$ such that,
for all $x, y$ in $E$, and $a$ in $A$,
$$\langle Tx,Ty\rangle=\varphi(\langle x,y\rangle)$$
and
$$T(xa)=T(x)\varphi(a).$$
\end{theorem}
\proof
We can regard $E$ and $F$ as full modules over $\langle E,E\rangle$ and $\langle F,F\rangle$, respectively.
In this case, as we mentioned above, $T$ is completely isometric and hence it preserves the triple products
$$T(z\langle x,y\rangle)=Tz\langle Tx,Ty\rangle, \forall x, y, z\in E.$$
Define $\varphi:\langle E,E\rangle\rightarrow\langle F,F\rangle$ by
$$\varphi(\sum\limits_{i=i}^n\alpha_i\langle x_i,y_i\rangle):=\sum\limits_{i=i}^n\alpha_i\langle Tx_i,Ty_i\rangle
,\ x_i, y_i\in E,\ \alpha_i\in\mathbb{C},\ i=1,\cdots,n.$$
Let $x_i, y_i$ and $z\in E,\ \alpha_i\in\mathbb{C},\ i=1,\cdots,n$. Then
$\sum\limits_{i=i}^n\alpha_i\langle x_i,y_i\rangle=0$ if and only if
$z(\sum\limits_{i=i}^n\alpha_i\langle x_i,y_i\rangle)=0$ for all $z$ if and only if
$T(z)(\sum\limits_{i=i}^n\alpha_i\langle Tx_i,Ty_i\rangle)=\sum\limits_{i=i}^n\alpha_iTz\langle Tx_i,Ty_i\rangle
=\sum\limits_{i=i}^n\alpha_iT(z\langle x_i,y_i\rangle)=T(z(\sum\limits_{i=i}^n\alpha_i\langle x_i,y_i\rangle))=0$
for all $z$
if and only if $\sum\limits_{i=i}^n\alpha_i\langle Tx_i,Ty_i\rangle=0$ since $T$ is injective,
$\sum\limits_{i=i}^n\alpha_i\langle x_i,y_i\rangle\in\langle E,E\rangle$ and
$\sum\limits_{i=i}^n\alpha_i\langle Tx_i,Ty_i\rangle\in\langle F,F\rangle$.
This shows that $\varphi$ is well-defined and injective.
From the defintion of $\varphi$, since $T$ is surjective, it is clear that $\varphi$ is a
surjective $*$-homomorphism and $T$ preserves the inner products w.r.t. $\varphi$.
By lemma \ref{lem:mp}, $T$ is a module map w.r.t $\varphi$.
\endproof

\begin{corollary}
Every surjective linear 2-isometry between two full Hilbert $C^*$-modules preserves the inner products and
is a module map with respect to some $*$-isomorphism of underlying $C^*$-algebras.
\end{corollary}

\bibliographystyle{amsplain}

\begin{thebibliography}{99}

\bibitem{Chu05}
C-H. Chu, M, Mackey, Isometries between ${\rm JB}^*$-triples.
\emph{Math. Z.} {\bf 251} (2005), no. 3, 615--633.

\bibitem{DupreGillette83}
M. J. Dupr\'{e} and R. M. Gillette, \emph{Banach bundles, Banach modules and automorphisms of $C^*$-algebras},
Research Notes in Mathematics 92, Pitman (1983).

\bibitem{Hamana99}
M. Hamana, Triple envelopes and Shilov boundaries of operator spaces,
\emph{Math. J. Toyama Univ.} {\bf 22} (1999), 77-93.

\bibitem{HsuWong11}
M. H. Hsu and N. C. Wong, Isometries embeddings of Banach bundles, to appear.

\bibitem{Isidro03}
J. M. Isidro, Holomorphic automorphisms of the unit balls of Hilbert
$C^*$-modules. \emph{Glasg. Math. J.} {\bf 45} (2003), no. 2, 249-262.

\bibitem{Jeang03}
J. S. Jeang and N. C. Wong, On the Banach-Stone Problem, {\em Studia Math.\ } {\bf 155} (2003), 95-105.

\bibitem{Jerison50}
M. Jerison,
The space of bounded maps into a Banach space, {\em Ann. of Math.\ } {\bf 52} (1950), 309-327.

\bibitem{Joita07}
M. Joi\c{t}a, A note about Hilbert modules over Fr\'{e}chet locally $C^*$-algebras,
\emph{Novi Sad J. Math.} {\bf 37} (2007), no. 1, 27-32.

\bibitem{Kadison51}
R. V. Kadison, Isometries of operator algebras, \emph{Ann. of Math.} {\bf 54} (1951), 325-338.

\bibitem{Kaup83}
W. Kaup,
A Riemann mapping theorem for bounded symmetric domains in complex Banach spaces,
\emph{Math. Z.} {\bf 183} (1983), 503-529.

\bibitem{Lance}
C. Lance, \emph{Hilbert $C^*$-modules}, London Mat. Soc. Lecture Notes Series, 210, cambridge University Press,
Cambridge, 1995.

\bibitem{Lau}
K. S. Lau,
A representation theorem for isometries of $C(X,E)$, {\em Pacific J. of Math.\ } {\bf 60} (1975), 229-233.

\bibitem{MS}
P. S. Muhly and B. Solel,
On the Morita equivalence of tensor algebras, \emph{Proc. London Math. Soc.} {\bf 81} (2000), 113-118.

\bibitem{NealRusso03}
M. Neal and B. Russo, Operator space characterizations of $C^*$-algebras and ternary rings,
\emph{Pac. J. Math.} {\bf 209} (2003), 339-364.

\bibitem{Pis03}
G. Pisier, \emph{Introduction to Operator Space Theory}, Cambridge University Press, 2003.

\bibitem{Solel01}
B. Solel,
Isometries of Hilbert $C^*$-modules, \emph{Trans. Amer. Math. Soc.} {\bf 553} (2001), 4637-4660.

\bibitem{Zettl83}
H. Zettl, A characterization of ternary rings of operators,
\emph{Adv. in Math.} {\bf }48 (1983), no. 2, 117-143.

\end{thebibliography}

\end{document}